# INVERSE NODAL PROBLEMS FOR DIRAC TYPE INTEGRO DIFFERENTIAL SYSTEM WITH A NONLOCAL BOUNDARY CONDITION

BAKI KESKIN

ABSTRACT. In this work, the Dirac-type integro differential system with one classical boundary condition and another nonlocal integral boundary condition is considered. We obtain the asymptotic formulae for the solutions, eigenvalues and nodal points. We also investigate the inverse nodal problem and prove that given a dense subset of nodal points uniquely determines the boundary condition parameter and the potential function of the considered differential system. We also provide an effective procedure for solving the inverse nodal problem.

## 1. Introduction

In recent years, boundary value problems with nonlocal conditions are one of the fastest developing current topics in mathematical physics. Nonlocal boundary conditions appear when we cannot measure data directly at the boundary. Problems of this type arise in various fields of mathematical physics, biology, biotechnology , mechanics, geophysics and other branches of natural sciences [2], [3], [4], [5], [6]. Such problems were studied firstly by Samarskii and Bitsadze. They formulated and investigated the non-local boundary problem for an elliptic equation [1].

Inverse spectral problems consist in recovering operators from their spectral characteristics. Results of the inverse problem for various nonlocal operators can be found in [7], [8], [9], [10], [43],[45]

In a classical inverse spectral problem, the potential and the coefficients of the operators are to be determined from spectral data (e.g., two sets of eigenvalues, or one set of eigenvalues and norming constants). This is called the inverse eigenvalue problem. An alternative is to take as data the zeros (nodes) of the eigenfunctions of the considered operators, which are just as experimentally observable as eigenvalues in some situations. This is generally referred to as the "Inverse Nodal Problem".

In 1988, Mclaughlin raised and solved the inverse nodal problem for Sturm-Liouville problems for the first time [35]. She showed that knowledge of a dense subset of zeros (nodes) of the eigenfunctions alone can determine the potential function of the Sturm-Liouville problem up to a constant. Hald and McLaughlin [25] provided some numerical schemes for the reconstruction of the potential for more general boundary conditions. In 1997, Yang suggested a constructive procedure

Current Address: Department of Mathematics, Faculty of Science, Sivas Cumhuriyet University 58140 Sivas, TURKEY
E-mail Address: bkeskin@cumhuriyet.edu.tr.

for reconstructing the potential and the boundary condition of the Sturm-Liouville problem from nodes of its eigenfunctions [41]. Inverse nodal problems have been addressed by various researchers in several papers for different operators with different kinds of boundary conditions ([16], [18], [21],[24] [37], [39], [40], [42], [44] and references therein). The inverse nodal problems for Dirac operators with various boundary conditions have been studied in [23], [38], [46], [47] and [48].

In recent years, integro-differential operators attracted much attention of mathematicians. Such operators have important applications in many fields of science (see monographs [17], [36] and references therein). Therefore, many researchers are currently working on inverse problems for these operators ([11], [12], [13], [14], [15], [19], [20],[22],[26] [27], [28], [30], [34] and [49]). The inverse nodal problem for Dirac type integro-differential operators with Robin boundary conditions was first studied by [29]. This operator with parameter dependent boundary conditions linearly and nonlinearly were studied by [31] and [32], respectively. In their studies, the authors considered $p(x)$ and $q(x)$, which are the components of the potential function $\Omega(x)$, as a special case such that $p(x) - q(x) = const$.

In this study, the perturbation of the Dirac differential system on a finite interval with one classical boundary condition and another nonlocal integral boundary condition by a Volterra type operator is considered. We also consider $p(x)$ and $q(x)$ as two independent functions and investigate a more general case. In this way, we have the opportunity to determine $p(x)$ and $q(x)$ separately. It is shown that the coefficients of the differential part of the operator and the boundary condition parameter can be determined by using a dense subset of the nodes. For this problem, we also give a constructive procedure for determining the coefficients as well as the useful asymptotics regarding the solution, eigenvalues and nodes.

## 2. Main Results

Consider the nonlocal boundary value problem of the Dirac type integro-differential system

$$BY' + \Omega(x)Y + \int_0^x M(x,t)Y\,dt = \lambda Y, \quad x \in (0,\pi), \tag{1}$$

with one classical boundary condition

$$y_1(0)\sin\theta + y_2(0)\cos\theta = 0 \tag{2}$$

and nonlocal integral boundary condition

$$y_1(\pi) = \int_0^\pi y_1(x)\omega(x)dx \tag{3}$$

where, $0 < \theta < \pi$ are real numbers, $\lambda$ is the spectral parameter,
$B = \begin{pmatrix} 0 & 1 \\ -1 & 0 \end{pmatrix}$, $\Omega(x) = \begin{pmatrix} p(x) & 0 \\ 0 & r(x) \end{pmatrix}$, $M(x,t) = \begin{pmatrix} M_{11}(x,t) & M_{12}(x,t) \\ M_{21}(x,t) & M_{22}(x,t) \end{pmatrix}$,
$Y = \begin{pmatrix} y_1 \\ y_2 \end{pmatrix}$, $p(x)$, $r(x)$, $\omega(x)$ and $M_{ij}(x,t)$, $(i,j=1,2)$ are real-valued and differentiable functions on segment $[0,\pi]$ f

Let $\varphi(x,\lambda) = (\varphi_1(x,\lambda), \varphi_2(x,\lambda))^T$ be the solution of (1) satisfying the initial condition $\varphi(0,\lambda) = (\cos\theta, -\sin\theta)^T$.

It is clear that $\varphi(x,\lambda)$ is an entire function of $\lambda$ satisfies the following Volterra integral equations:

$$
\begin{aligned}
\varphi_1(x,\lambda) &= \cos\theta\cos\lambda x + \sin\theta\sin\lambda x \\
&+ \int_0^x \sin\lambda(x-t)p(t)\varphi_1(t,\lambda)dt + \int_0^x \cos\lambda(x-t)r(t)\varphi_2(t,\lambda)d_\alpha t \\
&+ \int_0^x \int_0^t \sin\lambda(x-t)\{M_{11}(t,\xi)\varphi_1(\lambda,\xi) + M_{12}(t,\xi)\varphi_2(\lambda,\xi)\}d\xi dt \\
&+ \int_0^x \int_0^t \cos\lambda(x-t)\{M_{21}(t,\xi)\varphi_1(\lambda,\xi) + M_{22}(t,\xi)\varphi_2(\lambda,\xi)\}d\xi dt
\end{aligned}
\tag{4}
$$

$$
\begin{aligned}
\varphi_2(x,\lambda) &= \cos\theta\sin\lambda x - \sin\theta\cos\lambda x \\
&- \int_0^x \cos\lambda(x-t)p(t)\varphi_1(t,\lambda)dt + \int_0^x \sin\lambda(x-t)r(t)\varphi_2(t,\lambda)dt \\
&- \int_0^x \int_0^t \cos\lambda(x-t)\{M_{11}(t,\xi)\varphi_1(\lambda,\xi) + M_{12}(t,\xi)\varphi_2(\lambda,\xi)\}d\xi dt \\
&+ \int_0^x \int_0^t \sin\lambda(x-t)\{M_{21}(t,\xi)\varphi_1(\lambda,\xi) + M_{22}(t,\xi)\varphi_2(\lambda,\xi)\}d\xi dt
\end{aligned}
\tag{5}
$$

**Theorem 1.** *([28], the case: $\alpha = 1$) For $|\lambda| \to \infty$, the following asymptotic formulae are valid:*

$$
\begin{aligned}
\varphi_1(x,\lambda) &= \cos(\lambda x - \mu(x) - \theta) + \frac{1}{2\lambda}\upsilon(x)\cos(\lambda x - \mu(x) - \theta) \\
&- \frac{1}{2\lambda}\upsilon(0)\cos(\lambda x - \mu(x) + \theta) + \frac{1}{2\lambda}\sin(\lambda x - \mu(x) - \theta)\int_0^x \upsilon^2(t)dt \\
&- \frac{1}{2\lambda}K(x)\cos(\lambda x - \mu(x) - \theta) - \frac{1}{2\lambda}L(x)\sin(\lambda x - \mu(x) - \theta) \\
&+ o\left(\frac{1}{\lambda}\exp(|\tau|x)\right),
\end{aligned}
\tag{6}
$$

$$
\begin{aligned}
\varphi_2(x,\lambda) &= \sin(\lambda x - \mu(x) - \theta) - \frac{1}{2\lambda}\upsilon(x)\sin(\lambda x - \mu(x) - \theta) \\
&- \frac{1}{2\lambda}\upsilon(0)\sin(\lambda x - \mu(x) + \theta) - \frac{1}{2\lambda}\cos(\lambda x - \mu(x) - \theta)\int_0^x \upsilon^2(t)dt \\
&- \frac{1}{2\lambda}K(x)\sin(\lambda x - \mu(x) - \theta) + \frac{1}{2\lambda}L(x)\cos(\lambda x - \mu(x) - \theta) \\
&+ o\left(\frac{1}{\lambda}\exp(|\tau|x)\right),
\end{aligned}
\tag{7}
$$

*uniformly in* $x \in [0,\pi]$, *where,* $\mu(x) = \frac{1}{2}\int_0^x (p(t)+r(t))dt$, $\upsilon(x) = \frac{1}{2}(p(x)-r(x))$, $K(x) = \int_0^x (M_{11}(t,t) - M_{22}(t,t))dt$, $L(x) = \int_0^x (M_{12}(t,t) - M_{21}(t,t))dt$ *and* $\tau = \operatorname{Im}\lambda$.

The spectra of the boundary value problem (1)-(3) are the zero-sequences $\{\lambda_n\}_{n\in\mathbb{Z}}$ of the entire function

$$\Lambda(\lambda) := \varphi_1(\pi,\lambda) - \int_0^\pi \varphi_1(x,\lambda)\omega(x)dx = 0$$

and the function $\Lambda$ satisfies

$$\Lambda(\lambda) = \cos(\lambda\pi - \mu(\pi) - \theta) + \frac{1}{2\lambda}\upsilon(\pi)\cos(\lambda\pi - \mu(\pi) - \theta)$$
$$-\frac{1}{2\lambda}\upsilon(0)\cos(\lambda\pi - \mu(\pi) + \theta) + \frac{1}{2\lambda}\sin(\lambda\pi - \mu(\pi) - \theta)\int_0^\pi \upsilon^2(t)dt$$
(8)
$$-\frac{1}{2\lambda}K(\pi)\cos(\lambda\pi - \mu(\pi) - \theta) - \frac{1}{2\lambda}L(\pi)\sin(\lambda\pi - \mu(\pi) - \theta)$$
$$-\frac{1}{\lambda}\sin(\lambda\pi - \mu(\pi) - \theta)\omega(\pi) + o\left(\frac{1}{\lambda}\exp(|\tau|\pi)\right),$$

for sufficiently large $|\lambda|$. Since the eigenvalues of the problem (1)-(3) are the roots of $\Lambda(\lambda_n) = 0$, we can write the following equation for them:

$$\left(1 + \frac{1}{2\lambda_n}\upsilon(\pi) - \frac{1}{2\lambda_n}\upsilon(0)\cos 2\theta - \frac{1}{2\lambda_n}K(\pi)\right)\tan(\lambda_n\pi - \mu(\pi) - \frac{\pi}{2} - \theta) =$$
$$\frac{1}{2\lambda_n}\upsilon(0)\sin 2\theta + \frac{1}{2\lambda_n}\int_0^\pi \upsilon^2(t)dt - \frac{1}{2\lambda_n}L(\pi) - \frac{1}{\lambda_n}\omega(\pi) + o\left(\frac{1}{\lambda_n}\right)$$

from this expression, for sufficiently large $|n|$,

(9) $$\lambda_n = \left(n + \frac{1}{2}\right) + \frac{\theta + \mu(\pi)}{\pi}$$
$$+ \frac{1}{2n\pi}\left(\upsilon(0)\sin 2\theta + \int_0^\pi \upsilon^2(t)dt - L(\pi) - 2\omega(\pi)\right) + o\left(\frac{1}{n}\right), \; n \geq 1,$$

and similarly

(10) $$\lambda_n = \left(n - \frac{1}{2}\right) + \frac{\theta + \mu(\pi)}{\pi}$$
$$- \frac{1}{2n\pi}\left(\upsilon(0)\sin 2\theta + \int_0^\pi \upsilon^2(t)dt - L(\pi) - 2\omega(\pi)\right) + o\left(\frac{1}{n}\right), n \leq -1.$$

**Lemma 1.** *For sufficiently large $n$, $\varphi_1(x, \lambda_n)$ has exactly $n$ nodes $\{x_n^j : j = 0, 1, \cdots, n-1\}$ in the interval $(0, \pi)$: $0 < x_n^0 < x_n^1 < ... < x_n^{n-1} < \pi$. Moreover,*

$$x_n^j = \frac{(j + 1/2)\pi}{n} + \frac{\mu(x_n^j) + \theta}{n}$$
$$- \frac{(j + 1/2)\pi}{n}\frac{(\pi + 2\theta + 2\mu(\pi))}{2n\pi} - \left(\mu(x_n^j) + \theta\right)\frac{(\pi + 2\theta + 2\mu(\pi))}{2n^2\pi}$$
(11)
$$\frac{(j + 1/2)\pi}{n}\frac{1}{2n^2\pi}\left(\upsilon(0)\sin 2\theta + \int_0^\pi \upsilon^2(t)dt - L(\pi) - 2\omega(\pi)\right)$$
$$+ \frac{(j + 1/2)\pi}{n}\frac{(\pi + 2\theta + 2\mu(\pi))^2}{4n^2\pi^2} + \frac{1}{2n^2}\left(\upsilon(0)\sin 2\theta + \int_0^x \upsilon^2(t)dt - L(x)\right)$$
$$\frac{\mu(x_n^j) + \theta}{2n^3\pi}\left(\upsilon(0)\sin 2\theta + \int_0^\pi \upsilon^2(t)dt - L(\pi) - 2\omega(\pi)\right) -$$
$$+ \left(\mu(x_n^j) + \theta\right)\frac{(\pi + 2\theta + 2\mu(\pi))^2}{4n^3\pi^2} + O\left(\frac{1}{n^4}\right).$$

*uniformly with respect to $j \in \mathbb{Z}^+$.*

*Proof.* The first component $\varphi_1(x, \lambda_n)$ of the eigenfunction $\varphi(x, \lambda_n)$ has the following asymptotic formula:

$$\varphi_1(x, \lambda_n) = \cos(\lambda_n x - \mu(x) - \theta) + \frac{1}{2\lambda_n} v(x) \cos(\lambda_n x - \mu(x) - \theta)$$

(12)
$$-\frac{1}{2\lambda_n} v(0) \cos(\lambda_n x - \mu(x) + \theta) + \frac{1}{2\lambda_n} \sin(\lambda_n x - \mu(x) - \theta) \int_0^x v^2(t) dt$$

$$-\frac{1}{2\lambda_n} K(x) \cos(\lambda x - \mu(x) - \theta) - \frac{1}{2\lambda_n} L(x) \sin(\lambda_n x - \mu(x) - \theta)$$

$$+ o\left(\frac{1}{\lambda_n} \exp(|\tau| x)\right),$$

for $n \to \infty$ uniformly in $x$. Then for the nodal point $x_n^j$ of $\varphi_1(x, \lambda_n)$, from $\varphi_1(x_n^j, \lambda_n) = 0$, we obtain

$$\left(1 + \frac{1}{2\lambda_n} v(x_n^j) - \frac{1}{2\lambda_n} v(0) \cos 2\theta - \frac{1}{2\lambda_n} K(x_n^j)\right) \tan\left(\lambda_n x - \mu(x_n^j) - \theta - \frac{\pi}{2}\right) =$$

$$\frac{1}{2\lambda_n} v(0) \sin 2\theta + \frac{1}{2\lambda_n} \int_0^{x_n^j} v^2(t) dt - \frac{1}{2\lambda_n} L(x_n^j) + o\left(\frac{1}{\lambda_n}\right),$$

Taking into account Taylor's expansions, we get

$$\lambda_n x_n^j - \mu(x_n^j) - \theta - \frac{\pi}{2} = j\pi + \frac{1}{2\lambda_n}\left(v(0) \sin 2\theta + \int_0^{x_n^j} v^2(t) dt - L(x_n^j)\right) + o\left(\frac{1}{\lambda_n}\right).$$

It follows from the last equality

$$x_n^j = \frac{\left(j + \frac{1}{2}\right)\pi + \mu(x_n^j) + \theta}{\lambda_n} + \frac{1}{2\lambda_n^2}\left(v(0) \sin 2\theta + \int_0^{x_n^j} v^2(t) dt - L(x_n^j)\right) + o\left(\frac{1}{\lambda_n^2}\right).$$

The relation (13) is proven by using the asymptotic formula

$$\lambda_n^{-1} = \frac{1}{n}\left(1 - \frac{\pi + 2\theta + 2\mu(\pi)}{2n\pi} - \frac{\left(v(0) \sin 2\theta + \int_0^\pi v^2(t) dt - L(\pi) - 2\omega(\pi)\right)}{2n^2\pi}\right.$$

$$\left.\left(\frac{\pi + 2\theta + 2\mu(\pi)}{2n\pi}\right)^2 + o\left(\frac{1}{n^2}\right)\right)$$

as $n \to \infty$ uniformly in $j \in \mathbb{Z}^+$. □

**Theorem 2.** *Let $X$ be the set of nodal points. Fix $x \in (0, \pi)$. Let a sequence $\{x_n^j\} \subset X$ be chosen such that $x_n^j$ converges to $x$ as $n \to \infty$. Then the following limits are exist and finite and corresponding equalities hold:*

(13)
$$\lim_{|n| \to \infty} n\left(x_n^j - \frac{(j + 1/2)\pi}{n}\right) = \mu(x) + \theta - x\frac{\pi + 2\theta + 2\mu(\pi)}{2\pi} \triangleq f(x),$$

$$\lim_{|n|\to\infty} 2n^2\pi \left( x_n^j - \frac{(j+1/2)\pi}{n} + \frac{\mu(x_n^j)+\theta}{n} \right.$$

$$\left. + \frac{(j+1/2)\pi}{n} \frac{(\pi+2\theta+2\mu(\pi))}{2n\pi} \right)$$

(14)
$$= -(\mu(x)+\theta)(\pi+2\theta+2\mu(\pi))$$

$$-x\left(v(0)\sin 2\theta + \int_0^\pi v^2(t)dt - L(\pi) - 2\omega(\pi)\right)$$

$$+x\frac{(\pi+2\theta+2\mu(\pi))^2}{2\pi}$$

$$+\pi\left(v(0)\sin 2\theta + \int_0^x v^2(t)dt - L(x)\right) \triangleq g(x),$$

and

$$\lim_{|n|\to\infty} 2n^3\pi \left( x_n^j - \frac{(j+1/2)\pi}{n} + \frac{\mu(x_n^j)+\theta}{n} \right.$$

$$+ \frac{(j+1/2)\pi}{n}\frac{(\pi+2\theta+2\mu(\pi))}{2n\pi} + (\mu(x_n^j)+\theta)\frac{(\pi+2\theta+2\mu(\pi))}{2n^2\pi}$$

$$+ \frac{(j+1/2)\pi}{n}\frac{1}{2n^2\pi}\left(v(0)\sin 2\theta + \int_0^\pi v^2(t)dt - L(\pi) - 2\omega(\pi)\right)$$

(15)
$$- \frac{(j+1/2)\pi}{n}\frac{(\pi+2\theta+2\mu(\pi))^2}{4n^2\pi^2}$$

$$- \frac{1}{2n^2}\left(v(0)\sin 2\theta + \int_0^x v^2(t)dt - L(x)\right)\right)$$

$$= (\mu(x_n^j)+\theta)\left(v(0)\sin 2\theta + \int_0^\pi v^2(t)dt - L(\pi) - 2\omega(\pi)\right)$$

$$+ (\mu(x)+\theta)\frac{(\pi+2\theta+2\mu(\pi))^2}{2\pi} \triangleq h(x),$$

Therefore, proof of the following theorem is clear.

**Theorem 3.** *Let $\mu(\pi) = 0$. The given dense subset of nodal points $X$ uniquely determines the coefficients $\theta$ and $\omega(\pi)$ of the boundary conditions and if $L(x)$ is known, $X$ also uniquely determines the potential $\Omega(x)$ a.e. on $(0,\pi)$. Moreover, $\Omega(x), \omega(\pi)$ and $\theta$ can be via the following algorithm:*

*(1) For each $x \in (0,\pi)$ and $\alpha \in (0,1]$, choose a sequence $\{x_n^j\} \subset X$ such that $\lim_{|n|\to\infty} x_n^j = x$;*

*(2) Find the function $f(x)$ from (13) and calculate*

$$\theta = f(0)$$
$$2\mu'(x) = p(x) + r(x) = 2f'(x) + 2f(0) + \pi$$

*(3) Find the function $g(x)$ from (14) and calculate*

$$v(0) = \frac{g(0) - \theta\pi - 2\theta^2}{\pi \sin 2\theta}$$
$$\omega(\pi) = \frac{g(\pi) - 3\theta\pi - 4\theta^2 - \pi^2/2}{2\pi}$$

**(4)** If $L(x)$ is known then from (14) and (15) calculate
$$\begin{aligned} p(x) &= f'(x) + f(0) + \pi/2 + \rho(x) \\ r(x) &= f'(x) + f(0) + \pi/2 - \rho(x) \end{aligned}$$
where
$$\rho^2(x) = \frac{1}{\pi}\left(g'(x) + \left(f'(x) + f(0) + \frac{\pi}{2}\right)(\pi + 2\theta) + \frac{h(0)}{\theta} + L'(x)\right)$$

Thus, we have shown that we can reconstruct the potential function and obtain the coefficients of the boundary conditions using only dense subset of a nodal points. Our reconstruction formulae are also directly implies the uniqueness of this inverse problem.

**Conclusion:** In this work, we study inverse nodal problems for Dirac-type integro-differential systems with one classical boundary condition and another nonlocal integral boundary condition. We get useful asymptotics regarding the solution, eigenvalues, and nodes. And we present a constructive procedure to solve the inverse nodal problems. It makes sense for the theoretical integrity of the inverse nodal problem with nonlocal integral conditions.

**Conflict of Interest:** The author declares that there are no conflicts of interest regarding the publication of this paper

**Data Availability Statement:** Data sharing not applicable to this article as no datasets were generated or analysed during the current study

<em>8</em>       BAKI KESKIN